\documentclass[12pt]{amsproc}
\usepackage{amsmath,amssymb,amscd}

\begin{document}

\title[Lebesgue Decomposition of Functionals]{Lebesgue Decomposition of Functionals and Unique Preduals for Commutants Modulo Normed Ideals}
\author{Dan-Virgil Voiculescu}
\address{D.V. Voiculescu \\ Department of Mathematics \\ University of California at Berkeley \\ Berkeley, CA\ \ 94720-3840}
\thanks{Research supported in part by NSF Grant DMS-1301727.}
\date{}

\begin{abstract}
We prove an analogue of the Lebesgue decomposition for continuous functionals on the commutant modulo a reflexive normed ideal of an $n$-tuple of hermitian operators for which there are quasicentral approximate units relative to the normed ideal. Using results of Godefroy--Talagrand and Pfitzner we derive from this strong uniqueness of the predual of such a commutant modulo a normed ideal.
\end{abstract}

\maketitle

\section{Introduction}
\label{sec1}

On a separable infinite-dimensional complex Hilbert space ${\mathcal H}$ let $\tau = (T_1,\dots,T_n)$ be an $n$-tuple of bounded hermitian operators and let $({\mathcal J},|\cdot|_{\mathcal J})$ be a proper normed ideal of compact operators in the algebra ${\mathcal B}({\mathcal H})$ of bounded operators. The commutant of $\tau$ modulo ${\mathcal J}$ is the algebra ${\mathcal E}(\tau;{\mathcal J}) \subset {\mathcal B}({\mathcal H})$ of operators $X$ so that $[X,T_j] \in {\mathcal J}$, $1 \le j \le n$. With respect to the norm
\[
\||X\|| = \|X\| + \max_{1 \le j \le n} |[X,T_j]|_{\mathcal J}
\]
and the involution $X \to X^*$, ${\mathcal E}(\tau;{\mathcal J})$ is an involutive Banach algebra with isometric involution. For some of the Banach-space properties we will study $\||\cdot\||$ will be replaced by an equivalent norm.

Recently, motivated by operator theory problems, in \cite{21} we began studying certain algebras of this kind and then expanded the study in \cite{2}, \cite{22}, \cite{24} becoming aware that these algebras were natural objects of interest in connection with our old work (\cite{19}, \cite{20}) on normed ideal perturbation of Hilbert space operators. A basic insight of those papers was that our non-commutative Weyl--von~Neumann type theorem (\cite{18}) could be adapted to deal with proper normed ideals (\cite{19}) and that this was the starting point for a new approach also for other perturbation theory questions (see the survey paper \cite{20}). Another source for renewed interest was the paper \cite{3} of A.~Connes where our perturbation results are used in dealing with a non-commutative geometry question.

Assuming $k_{\mathcal J}(\tau) = 0$, a condition which amounts to the existence of $0 \le A_m \le I$, $A_m \uparrow I$ finite rank operators so that
\[
\lim_{m \to \infty} \max_{1 \le j \le n} |[A_m,T_j]|_{\mathcal J} = 0,
\]
often accompanied by the assumption that the finite rank operators are dense in ${\mathcal J}$ and in its dual ${\mathcal J}^d$, we exhibited many similarities between the pair $({\mathcal K}(\tau;{\mathcal J}),{\mathcal E}(\tau;{\mathcal J}))$, where ${\mathcal K}(\tau;{\mathcal J})$ denotes the compact operators in ${\mathcal E}(\tau;{\mathcal J})$, and the pair $({\mathcal K}({\mathcal H}),{\mathcal B}({\mathcal H}))$, where ${\mathcal K}({\mathcal H})$ denotes the compact operators in ${\mathcal B}({\mathcal H})$. In particular, ${\mathcal E}(\tau;{\mathcal J})/{\mathcal K}(\tau;{\mathcal J})$ like the Calkin algebra ${\mathcal B}({\mathcal H})/{\mathcal K}({\mathcal H})$ is a $C^*$-algebra and we also proved the analogues of ${\mathcal B}({\mathcal H})$ being the second dual of ${\mathcal K}({\mathcal H})$ and the multiplier algebra of ${\mathcal K}({\mathcal H})$. This leads immediately to the two questions we address here: the analogues of the uniqueness of the predual of ${\mathcal B}({\mathcal H})$ and the decomposition of a functional as a sum of an ultraweakly continuous and of a singular functional (the particular case of the von~Neumann algebra ${\mathcal B}({\mathcal H})$ of theorems of S.~Sakai \cite{14} and M.~Takesaki \cite{16}). Thus ${\mathcal E}(\tau;{\mathcal J})$ with a conveniently modified norm joins a growing list of Banach spaces and algebras for which these kind of properties have been established (starting with the work of Grothendieck \cite{9}, followed by \cite{14}, \cite{16}, \cite{1}, \cite{7}, \cite{6}, \cite{13}, \cite{5}, \cite{12}, \cite{17}, \cite{4}, \cite{11}).

In addition to the introduction, there are three more sections. In Section~\ref{sec2} we collect the necessary preliminaries. Then in Section~\ref{sec3} we establish the Lebesgue decomposition of functionals on ${\mathcal E}(\tau;{\mathcal J})$ and in Section~\ref{sec4} we prove the uniqueness of the predual. The uniqueness of the predual is easier to establish, being able to use condition $(X)$ of Godefroy and Talagrand (\cite{7}) and a theorem of Pfitzner (\cite{12}).

\section{Preliminaries}
\label{sec2}

\subsection{Normed ideals}
\label{subsec2.1}

By a {\em normed ideal} we mean what usually goes under the longer name of a {\em symmetrically normed ideal} (\cite{10}, \cite{15}). These are ideals $0 \ne {\mathcal J} \subset {\mathcal B}({\mathcal H})$, that are contained in the ideal ${\mathcal K}({\mathcal H})$ of compact operators and that are endowed with a norm $|\cdot|_{\mathcal J}$ satisfying certain conditions, like $|AXB|_{\mathcal J} \le \|A\|\,|X|_{\mathcal J}\|B\|$ if $A,B \in {\mathcal B}({\mathcal H})$ and $X \in {\mathcal J}$ and with respect to which ${\mathcal J}$ is a Banach space (for more see \cite{10}, \cite{15}). We shall call ${\mathcal J}$ a {\em proper normed ideal} if ${\mathcal J} \ne {\mathcal K}({\mathcal H})$. The norm can be written $|T|_{\mathcal J} = |T|_{\Phi} = \Phi(s_1(T),s_2(T),\dots)$ where $\Phi$ is a {\em norming function} (see \S 3 in \cite{10}) and $s_1(T) \ge s_2(T) \ge \dots$ are the eigenvalues of the compact operator $(T^*T)^{1/2}$. Like in \cite{10}, given a norming function $\Phi$, we denote by ${\mathcal G}_{\Phi}$ the ideal of all compact operators $T$ such that $|T|_{\Phi} < \infty$ and by ${\mathcal G}_{\Phi}^{(0)}$ the closure in ${\mathcal G}_{\Phi}$ of the ideal ${\mathcal R}({\mathcal H})$ of finite-rank operators. If ${\mathcal G}_{\Phi} = {\mathcal G}_{\Phi}^{(0)}$, the norming function $\Phi$ is called ``mononorming'' (\cite{10}). By $({\mathcal C}_p,|\ |_p)$, $1 \le p < \infty$ we denote the Schatten--von~Neumann $p$-class. The norming function of ${\mathcal C}_p$ is mononorming if $1 \le p < \infty$. If $\Phi$ is a norming function, the dual of the Banach space $({\mathcal G}_{\Phi}^{(0)},|\ |_{\Phi})$ is $({\mathcal G}_{\Phi^*},|\ |_{\Phi^*})$ under the duality $(X,Y) \to Tr(XY)$ for $(X,Y) \in {\mathcal G}_{\Phi}^{(0)} \times {\mathcal G}_{\Phi^*}$ (we leave out the case ${\mathcal G}_{\Phi}^{(0)} = {\mathcal C}_1$ where the dual is ${\mathcal B}({\mathcal H})$). Here $\Phi^*$ is the conjugate norming function to $\Phi$ (see \cite{10}). In case both $\Phi$ and $\Phi^*$ are mononorming, ${\mathcal G}_{\Phi}$ and ${\mathcal G}_{\Phi^*}$ are each other's dual and are reflexive Banach spaces (Thm.~2.2, Ch.~III in \cite{10}).

\subsection{The condition $k_{\Phi}(\tau) = 0$}
\label{subsec2.2}

If $\Phi$ is a norming function and $\kappa = (K_1,\dots,K_n) \in ({\mathcal K}({\mathcal H}))^n$ by $|\kappa|_{\Phi}$ we shall denote the $\max_{1 \le j \le n} |K_j|_{\Phi}$ and if $\tau = (T_1,\dots,T_m) \in ({\mathcal B}({\mathcal H}))^n$ and $X \in {\mathcal B}({\mathcal H})$, by $[\tau,X]$ we denote the $n$-tuple $([T_1,X],\dots,[T_n,X])$. The set ${\mathcal R}^+_1({\mathcal H})$ of finite rank operators $A$ such that $0 \le A \le I$ is a directed set with respect to $\le$ and if $\tau \in ({\mathcal B}({\mathcal H}))^n$ then
\[
k_{\Phi}(\tau) = \liminf_{A \in {\mathcal R}^+_1({\mathcal H})} |[\tau,A]|_{\Phi}
\]
(see \cite{20} for a survey of our work on this invariant). The condition $k_{\Phi}(\tau) = 0$ is equivalent to the asymptotic commutation condition we mentioned in the introduction: there exist $A_j \in {\mathcal R}_1^+({\mathcal H})$ so that $A_j \uparrow I$ as $j \to \infty$ and $\lim_{j \to \infty}|[\tau,A_j]|_{\Phi} = 0$.

\subsection{The algebras ${\mathcal E}(\tau;{\mathcal J})$, ${\mathcal K}(\tau;{\mathcal J})$, ${\mathcal E}/{\mathcal K}(\tau;{\mathcal J})$}
\label{subsec2.3}

Given a normed ideal $({\mathcal J},|\ |_{\mathcal J})$ and $\tau \in ({\mathcal B}({\mathcal H}))^n$, $\tau = \tau^*$ an $n$-tuple of hermitian operators, we consider the algebras (see \cite{22})
\[
\begin{aligned}
{\mathcal E}(\tau;{\mathcal J}) &= \{S \in {\mathcal B}({\mathcal H}) \mid [T_j,S] \in {\mathcal J},\ 1 \le j \le n\}, \\
{\mathcal K}(\tau;{\mathcal J}) &= {\mathcal E}(\tau;{\mathcal J}) \cap {\mathcal K}({\mathcal H}), \\
{\mathcal E}/{\mathcal K}(\tau;{\mathcal J}) &= {\mathcal E}(\tau;{\mathcal J})/{\mathcal K}(\tau;{\mathcal J}).
\end{aligned}
\]
With the norm
\[
\||S\|| = \|S\| + |[\tau,S]|_{\mathcal J}
\]
${\mathcal E}(\tau;{\mathcal J})$ becomes a Banach algebra with an isometric involution and ${\mathcal K}(\tau;{\mathcal J})$ is a closed ideal, while ${\mathcal E}/{\mathcal K}(\tau;{\mathcal J})$ identifies algebraically with a $*$-subalgebra of the Calkin algebra ${\mathcal B}/{\mathcal K}({\mathcal H}) = {\mathcal B}({\mathcal H})/{\mathcal K}({\mathcal H})$.

We showed (Prop.~2.5 of \cite{22}) that under the assumptions that ${\mathcal J} = {\mathcal G}_{\Phi}^{(0)}$ and $k_{\Phi}(\tau) = 0$, where $\Phi$ is the norming function so that $|\ |_{\mathcal J} = |\ |_{\Phi}$ on ${\mathcal J}$, we have that ${\mathcal E}/{\mathcal K}(\tau;{\mathcal J})$ actually identifies isometrically with the corresponding $*$-subalgebra of the Calkin algebra, so that this $*$-subalgebra is actually a $C^*$-subalgebra of ${\mathcal B}/{\mathcal K}$ and also ${\mathcal E}/{\mathcal K}(\tau;{\mathcal J})$ with the quotient norm is a $C^*$-algebra. The condition ${\mathcal J} = {\mathcal G}_{\Phi}^{(0)}$ means that ${\mathcal R}({\mathcal H})$ is dense in ${\mathcal J}$ and this is automatic when $\Phi$ is mononorming. We shall also use the notation ${\mathcal E}(\tau;\Phi)$, ${\mathcal K}(\tau;\Phi)$, ${\mathcal E}/{\mathcal K}(\tau;\Phi)$ for ${\mathcal E}(\tau;{\mathcal J})$, ${\mathcal K}(\tau;{\mathcal J})$, ${\mathcal E}/{\mathcal K}(\tau;{\mathcal J})$ when $\Phi$ is mononorming.

Besides the norm $\||\ \||$ on ${\mathcal E}(\tau;{\mathcal J})$ it is also useful to introduce the equivalent Banach space norm
\[
\||S\||_M = \max(\|S\|,|[\tau,S]|_{\Phi}).
\]
Note that a similar max-norm is used for Lipschitz-algebras (\cite{25}).

\subsection{Dualities}
\label{sec2.4}

{\em The duality results in this subsection are all under the assumption that $\Phi$ and $\Phi^*$ are mononorming and that $k_{\Phi}(\tau) = 0$.}

Let ${\mathcal N}$ be the subspace of ${\mathcal C}_1 \times ({\mathcal G}_{\Phi^*})^n$ consisting of those vectors $(x,(y_j)_{1 \le j \le n})$ such that $x = \sum_{1 \le j \le n} [T_j,y_j]$. Consider the three Banach spaces
\[
{\mathcal K}(\tau;\Phi),({\mathcal C}_1 \times ({\mathcal G}_{\Phi^*})^n)/{\mathcal N},{\mathcal E}(\tau;\Phi).
\]
We proved in \cite{22} that the dual of the first identifies with the second, while the dual of the second identifies with the third, so that in particular ${\mathcal E}(\tau;\Phi)$ identifies with the bidual of ${\mathcal K}(\tau;\Phi)$. The duality pairings are those which arise from identifying ${\mathcal K}(\tau;\Phi)$ and ${\mathcal E}(\tau;\Phi)$ with subspaces of ${\mathcal K}({\mathcal H}) \times ({\mathcal G}_{\Phi})^n$ and ${\mathcal B}({\mathcal H}) \times ({\mathcal G}_{\Phi})^n$, respectively, and the fact that ${\mathcal C}_1 \times ({\mathcal G}_{\Phi^*})^n$ is the dual of ${\mathcal K}({\mathcal H}) \times ({\mathcal G}_{\Phi})^n$ and the predual of ${\mathcal B}({\mathcal H}) \times ({\mathcal G}_{\Phi})^n$ via the trace-pairings. The isometric embeddings of ${\mathcal K}(\tau;\Phi)$ and ${\mathcal E}(\tau;\Phi)$ take an operator $S$ to $(S,([T_j,S])_{1 \le j \le n})$ viewed as an element of ${\mathcal K}({\mathcal H}) \times ({\mathcal G}_{\Phi})^n$ or ${\mathcal B}({\mathcal H}) \times ({\mathcal G}_{\Phi})^n$, endowed with the norm
\[
\|(S,(X_j)_{1 \le j \le n})\| = \|S\| + \max_{1 \le j \le n} |X_j|_{\Phi}.
\]
The norm on ${\mathcal C}_1 \times ({\mathcal G}_{\Phi^*})^n$ is here
\[
\|(X,(Y_j)_{1 \le j \le n})\| = \max\left(|X|_1,\sum_{1 \le j \le n}|Y_j|_{\Phi^*}\right).
\]
(The fact that the identification of ${\mathcal E}(\tau;\Phi)$ with the dual of $({\mathcal C}_1 \times ({\mathcal G}_{\Phi^*})^n)/{\mathcal N}$ is actually an isometric identification is obvious from the embedding once the dualities have been established.)

In the present paper we shall use the notation ${\mathcal E}_*(\tau;\Phi)$ for $({\mathcal C}_1 \times ({\mathcal G}_{\Phi^*})^n)/{\mathcal N}$, while ${\mathcal E}^*(\tau;\Phi)$ will denote the dual of ${\mathcal E}(\tau;\Phi)$. Then ${\mathcal E}_*(\tau;\Phi)$ identifies isometrically with the subspace of ${\mathcal E}^*(\tau;\Phi)$ consisting of functionals $\varphi$ on ${\mathcal E}(\tau;\Phi)$ such that
\[
\varphi(S) = Tr\left(SX + \sum_{1 \le j \le n} Y_j[T_j,S]\right)
\]
where $X \in {\mathcal C}_1$ and $Y_j \in {\mathcal G}_{\Phi^*}$, $1 \le j \le n$. The fact that ${\mathcal E}_*(\tau;\Phi)$ is the dual of ${\mathcal K}(\tau;\Phi)$ means that the restriction of a functional $\varphi \in {\mathcal E}^*(\tau;\Phi)$ to ${\mathcal K}(\tau;\Phi)$ yields a functional in ${\mathcal E}_*(\tau;\Phi)$.

If we use on ${\mathcal E}(\tau;\Phi)$ and ${\mathcal K}(\tau;\Phi)$ the equivalent norm
\[
\||S\||_M = \max(\|S\|,|[T_1,S]|_{\Phi},\dots,|[T_n,S]|_{\Phi})
\]
then to keep the identifications isometric, the norm on ${\mathcal E}_*(\tau;\Phi)$ should be replaced by the quotient norm on $({\mathcal C}_1 \times ({\mathcal G}_{\Phi^*})^n)/{\mathcal N}$ where ${\mathcal C}_1 \times ({\mathcal G}_{\Phi^*})^n$ is now endowed with the norm
\[
\|(x,(y_j)_{1 \le j \le n})\|_M = |x|_1 + \sum_{1 \le j \le n} |Y_j|_{\Phi^*}.
\]

\section{The Lebesgue decomposition}
\label{sec3}

\subsection{Weak and ultraweak functionals}
\label{sec3.1}

We begin the discussion of functionals in a quite general setting, but we will have to impose various restrictive conditions before long.

Let $\tau \in ({\mathcal B}({\mathcal H}))^n$ and let $({\mathcal J},|\ |_{\mathcal J})$ be a normed ideal. A functional $\varphi: {\mathcal E}(\tau;{\mathcal J}) \to {\mathbb C}$ will be called a {\em weak functional} if there are finite rank operators $A,B_1,\dots,B_n \in {\mathcal R}({\mathcal H})$ so that
\[
\varphi(S) = Tr\left(SA + \sum_{1 \le j \le n} [T_j,S]B_j\right).
\]
Note that since
\[
\varphi(S) = Tr\left( S\left( A - \sum_{1 \le j \le n} [T_j,B_j]\right)\right)
\]
we can also do without the $B_j$'s, that is $\varphi(S) = Tr(SA)$. Assuming ${\mathcal J} = {\mathcal G}_{\Phi}^{(0)}$ for some norming function $\Phi$ then $\varphi$ will be called an {\em ultraweak fuctional} if there are $A \in {\mathcal C}_1,B_1,\dots,B_n \in {\mathcal G}_{\Phi^*}$ so that
\[
\varphi(S) = Tr\left( SA + \sum_{1 \le j \le n} [S,T_j]B_j\right).
\]
Using the isometric embeddings of ${\mathcal E}(\tau;{\mathcal J})$ into ${\mathcal B}({\mathcal H}) \times ({\mathcal G}_{\Phi}^{(0)})^n$ for the norms $\||\cdot\||$ and $\||\cdot\||_M$, respectively, we easily see that
\[
\|\varphi\| \le \max\left( |A|_1,\sum_{1 \le j \le n} |B_j|_{\Phi^*}\right)
\]
and
\[
\|\varphi\|_M \le |A|_1 + \sum_{1 \le j \le n} |B_j|_{\Phi^*}.
\]
There is also another natural class of functionals $\varphi: {\mathcal E}(\tau;{\mathcal J}) \to {\mathbb C}$, where ${\mathcal G}_{\Phi}^{(0)} \subset {\mathcal J} \subset {\mathcal G}_{\Phi}$ and we have
\[
\varphi(S) = Tr\left( SA + \sum_{1 \le j \le n} [S,T_j]B_j\right)
\]
where $A \in {\mathcal C}_1$, $B_j \in {\mathcal G}_{\Phi^*}^{(0)}$ which in the absence of a better name we may call {\em intermediate weak functionals}. Clearly the same inequalities for $\|\varphi\|$ and $\|\varphi\|_M$ hold. The weak functionals are a dense subset in the norm-topology of the set of intermediate weak functionals. In case $\Phi$ and $\Phi^*$ are mononorming, the ultraweak and intermediate weak functionals coincide.

\subsection{Singular functionals}
\label{sec3.2}

In general a functional $\varphi: {\mathcal E}(\tau;{\mathcal J}) \to {\mathbb C}$ will be called a {\em singular functional} if there is a continuous functional $\psi: {\mathcal B}({\mathcal H}) \to {\mathbb C}$ so that $\psi({\mathcal K}({\mathcal H})) = 0$ and $\varphi(S) = \psi(S)$.

The decomposition of functionals on ${\mathcal B}({\mathcal H})$ then easily gives the following result.

\bigskip
\noindent
{\bf 3.1. Proposition.} {\em 
Let $\varphi: {\mathcal E}(\tau;{\mathcal G}_{\Phi}^{(0)}) \to {\mathbb C}$ be a continuous functional. Then there is a continuous functional $\psi: {\mathcal B}({\mathcal H}) \to {\mathbb C}$ so that $\psi({\mathcal K}({\mathcal H})) = 0$ and there are $A \in {\mathcal C}_1$ and $B_j \in {\mathcal G}_{\Phi^*}$, $1 \le j \le n$ so that
\[
\varphi(S) = \psi(S) + Tr\left( SA + \sum_{1 \le j \le n} [S,T_j]B_j\right)
\]
and
\[
\|\varphi\|_M = \|\psi\| + |A|_1 + \sum_{1 \le j \le n} |B_j|_{\Phi^*}.
\]
Similarly, $A$ and $B_j$, $1 \le j \le n$, can be chosen so that
\[
\|\varphi\| = \max\left( \|\psi\| + |A|_1, \sum_{1 \le j \le n} |B_j|_{\Phi^*}\right).
\]
}

\bigskip
\noindent
{\bf {\em Proof.}} We embed isometrically ${\mathcal E}(\tau;{\mathcal G}_{\Phi}^{(0)})$ into ${\mathcal B}({\mathcal H}) \times ({\mathcal G}_{\Phi}^{(0)})^n$ by sending $S$ to $(S,([S,T_j])_{1 \le j \le n})$ and endowing ${\mathcal B}({\mathcal H}) \times ({\mathcal G}_{\Phi}^{(0)})^n$ with the appropriate norm in case we deal with $\||\cdot\||$ or $\||\ \||_M$ on ${\mathcal E}(\tau;{\mathcal G}_{\Phi}^{(0)})$. By Hahn--Banach $\varphi$ can be extended to ${\mathcal B}({\mathcal H}) \times ({\mathcal G}_{\Phi}^{(0)})^n$ with preservation of its norm and we get a functional ${\tilde \psi}: {\mathcal B}({\mathcal H}) \to {\mathbb C}$ and $B_j \in {\mathcal G}_{\Phi^*}$, $1 \le j \le n$, so that
\[
\varphi(S) = {\tilde \psi}(S) + Tr\left( \sum_{1 \le j \le n} [T_j,S]B_j\right)
\]
and
\[
\|\varphi\| = \max\left( \|{\tilde \psi}\|,\sum_{1 \le j \le n} |B_j|_{\Phi^*}\right)
\]
or
\[
\|\varphi\|_M = \|{\tilde \psi}\| + \sum_{1 \le j \le n} |B_j|_{\Psi^*}
\]
depending on the choice of norm on ${\mathcal E}(\tau;{\mathcal G}_{\Phi}^{(0)})$. The decomposition \cite{16} of ${\tilde \psi}(S) = \psi(S) + Tr(AS)$ where $\psi: {\mathcal B}({\mathcal H}) \to {\mathbb C}$ is so that $\psi({\mathcal K}({\mathcal H})) = 0$, $A \in {\mathcal C}_1$ and $\|{\tilde \psi}\| = \|\psi\| + |A|_1$ then concludes the proof.\qed

\bigskip
\noindent
{\bf 3.1. Corollary.} {\em 
If $\varphi: {\mathcal E}(\tau;{\mathcal G}_{\Phi}^{(0)}) \to {\mathbb C}$ is a continuous functional, then there is a singular functional $\varphi_1$ and an ultraweak functional $\varphi_2$ so that $\varphi = \varphi_1 + \varphi_2$ and $\|\varphi\|_M = \|\varphi_1\|_M + \|\varphi_2\|_M$.
}

\bigskip
\noindent
{\bf {\em Proof.}} Clearly, with the $\psi,A$ and $B_j - s$ which we found in Proposition~3.1 for the norm $\||\ \||_M$ putting $\varphi_1(S) = \psi(S)$ and
\[
\varphi_2(S) = Tr\left( SA + \sum_{1 \le j \le n} [S,T_j]B_j\right)
\]
we will have $\varphi = \varphi_1 + \varphi_2$ and since $\|\varphi_1\|_M \le \|\psi\|$,
\[
\|\varphi_2\|_M \le |A|_1 + \sum_{1 \le j \le n} |B_j|_{\Phi^*}
\]
we will have $\|\varphi\|_M \ge \|\varphi_1\|_M + \|\varphi_2\|_M$ and this must be an equality, because of the triangle inequality.\qed

\bigskip
{\em Under the additional assumption that $k_{\Phi}(\tau) = 0$, we will prove uniqueness of the decomposition, which makes the result substantially stronger. We will call this unique decomposition the Lebesgue decomposition of $\varphi$.}

\bigskip
\noindent
{\bf 3.2. Proposition.} {\em 
Assume $k_{\Phi}(\tau) = 0$. Then, if $\varphi: {\mathcal E}(\tau;{\mathcal G}_{\Phi}^{(0)}) \to {\mathbb C}$ is a continuous functional, there is a unique decomposition $\varphi = \varphi_a + \varphi_s$ where $\varphi_a$ is an ultraweak functional and $\varphi_s$ is a singular functional. Moreover, we have $\|\varphi\|_M = \|\varphi_a\|_M + \|\varphi_s\|_M$.
}

\bigskip
\noindent
{\bf {\em Proof.}} In view of Corollary~3.1 which established the existence of a decomposition which satisfies the relation among norms, it clearly suffices to prove the uniqueness, that is a continuous functional $\varphi$ that is at the same time singular and ultraweak must be $= 0$. Since $k_{\Phi}(\tau) = 0$, there are $A_k \in {\mathcal R}^+_1({\mathcal H})$ so that $A_k \uparrow I$ and $|[A_k,\tau]|_{\Phi} \to 0$ as $k \to \infty$. If $S \in {\mathcal E}(\tau;\Phi)$ we have $\varphi(A_kS) = 0$ since $\varphi$ is singular. On the other hand, since $\varphi$ is ultraweak, there are $A \in {\mathcal C}_1$, $B_j \in {\mathcal G}_{\Phi^*}$ so that
\[
\varphi(S) = Tr\left( AS + \sum_{1 \le j \le n} [S,T_j]B_j\right).
\]
We shall prove that $\lim_{k \to \infty} \varphi(A_kS) = \varphi(S)$. Since $A_k \uparrow I$, we have $\lim_{k \to \infty} |A - AA_k|_1 = 0$ and hence $\lim_{k \to \infty} Tr\ AA_kS = Tr\ AS$. Further we have $Tr[A_kS,T_j]B_j = Tr[A_k,T_j]SB_j + Tr\ A_k[S,T_j]B_j$ and $|[A_k,T_j]SB_j|_1 \le |[A_k,\tau]|_{\Phi} |SB_j|_{\Phi} \le \|S\| |[A_k,\tau]|_{\Phi} |B_j|_{\Phi^*}$ so that
\[
\lim_{k \to \infty} |[A_k,T_j]SB_j|_1 = 0
\]
and also
\[
|A_k[S,T_j]B_j - [S,T_j]B_j|_1 \to 0
\]
as $k \to \infty$ since $[S,T_j]B_j \in {\mathcal C}_1$. This gives $\lim_{k \to \infty} \varphi(A_kS) = \varphi(S)$ and hence $\varphi(S) = 0$.\qed

\bigskip
The limit underlying the proof of the preceding proposition upon closer examination can be used to describe $\varphi_a$ when $\varphi$ is given.

\bigskip
\noindent
{\bf 3.3. Proposition.} {\em 
Assume $k_{\Phi}(\tau) = 0$ and let $\varphi: {\mathcal E}(\tau;{\mathcal G}_{\Phi}^{(0)}) \to {\mathbb C}$ with Lebesgue decomposition $\varphi = \varphi_a + \varphi_s$. Let further $A_k \in {\mathcal R}^+_1({\mathcal H})$, $k \in {\mathbb N}$ be a sequence so that $A_k \uparrow I$ and $|[A_k,\tau]|_{\Phi} \to 0$ as $k \to \infty$. Then we have
\[
\varphi_a(S) = \lim_{k \to \infty} \varphi(A_kS)
\]
for all $S \in {\mathcal E}(\tau;{\mathcal G}_{\Phi}^{(0)})$. Moreover, if $\Phi^*$ is mononorming and we define $\varphi_k$ by $\varphi_k(S) = \varphi(A_kS)$ then
\[
\lim_{k \to \infty} \|\varphi_a - \varphi_k\| = 0.
\]
}

\bigskip
\noindent
{\bf {\em Proof.}} Since $\varphi_s(A_kS) = 0$ it suffices to prove the assertions in the case $\varphi = \varphi_a$. Thus 
\[
\varphi(S) = Tr\left(AS + \sum_{1 \le j \le n} [S,T_j]B_j\right)
\]
where $A \in {\mathcal C}_1$ and $B_j \in {\mathcal G}_{\Phi^*}$, $1 \le j \le n$. We have
\[
\begin{aligned}
|\varphi_a(S) - \varphi(A_kS)| &\le |Tr(A-AA_k)S| + \sum_{1 \le j \le n} |Tr((A_k[S,T_j] - [S,T_j])B_j)| \\
&\quad + \sum_{1 \le j \le n} |Tr([A_k,T_j]SB_j)| \\
&\le |A - AA_k|_1\|S\| + \sum_{1 \le j \le n} |(I-A_k)[S,T_j]|_{\Phi}|B_j|_{\Phi^*} \\
&\quad + \sum_{1 \le j \le n} |[A_k,\tau]|_{\Phi} |B_j|_{\Phi^*} \|S\|.
\end{aligned}
\]
Clearly $|A-AA_k|_1 \to 0$, $|(I-A_k)[S,T_j]|_{\Phi} \to 0$, $|[A_k,\tau]|_{\Phi} \to 0$ as $k \to \infty$, where we used that $[S,T_j] \in {\mathcal G}_{\Phi}^{(0)}$. To prove the second assertion, remark that if $\Phi^*$ is mononorming then
\[
\begin{aligned}
|Tr(A_k-I)[S,T_j]B_j| &= |Tr[S,T_j]B_j(A_k-I)| \\
&\le |[S,\tau]|_{\Phi} |B_j(A_k-I)|_{\Phi^*} \\
&\le \||S\|| |B_j(A_k-I)|_{\Phi^*}
\end{aligned}
\]
and $|B_j(A_k-I)|_{\Phi^*} \to 0$ as $k \to \infty$ because $\Phi^*$ is mononorming. Thus the earlier estimate gives now
\[
\begin{aligned}
&|\varphi_a(S)-\varphi(A_kS)| \\
&\quad \le \||S\||\left( |A-AA_k|_1 + \sum_{1 \le j \le n} (|B_j(A_k-I)|_{\Phi^*} + |[\tau,A_k]|_{\Phi}|B_j|_{\Phi^*})\right)
\end{aligned}
\]
and the quantity multiplying $\||S\||$ converges to $0$ as $k \to \infty$, which gives the desired result.\qed

\bigskip
\noindent
{\bf 3.1. Remark.} Assuming $k_{\Phi}(\tau) = 0$ we know that ${\mathcal E}/{\mathcal K}(\tau;{\mathcal G}_{\Phi}^{(0)})$ identifies isometrically with a $C^*$-subalgebra of the Calkin algebra. It is easily seen that the singular functionals in this case are precisely the functionals on ${\mathcal E}(\tau;{\mathcal G}_{\Phi}^{(0)})$ which arise from functionals on the $C^*$-algebra ${\mathcal E}/{\mathcal K}(\tau;{\mathcal G}_{\Phi}^{(0)})$ and that the $\|\ \|$ and $\|\ \|_M$ norms of these functionals coincide with the norm of the functional on the $C^*$-algebra ${\mathcal E}/{\mathcal K}(\tau;{\mathcal G}_{\Phi}^{(0)})$.

\section{Unique preduals}
\label{sec4}

\subsection{Preduals}
\label{subsec4.1}

{\em Throughout Section $4$ we will assume that $\Phi$ and $\Phi^*$ are mononorming and that $k_{\Phi}(\tau) = 0$.} These are the same assumptions like in $2.4$, where the preliminaries about dualities were presented. Thus ${\mathcal K}(\tau;\Phi)$ has a dual ${\mathcal E}_*(\tau;\Phi)$ and the dual of ${\mathcal E}_*(\tau;\Phi)$ is ${\mathcal E}(\tau;\Phi)$, these dualities being isometric with respect to the norms $\||\cdot\||$, $\|\cdot\|$, $\||\cdot\||$ or with respect to the equivalent norms $\||\ \||_M$, $\|\ \|_M$, $\||\ \||_M$. In view of the results about the Lebesgue decomposition, we have a slight preference for the second set of norms. In view of the description of the duality between ${\mathcal E}_*(\tau;\Phi)$ and ${\mathcal E}(\tau;\Phi)$ (\cite{22} Prop.~$4.5$) it is {\em immediate that ${\mathcal E}_*(\tau;\Phi)$ viewed as a subspace of ${\mathcal E}^*(\tau;\Phi)$ identifies with the ultraweak functionals on ${\mathcal E}(\tau;\Phi)$}.

\subsection{The $L$-embedded subspace property}
\label{subsec4.2}

We recall that a Banach-space ${\mathcal X}$ is $L$-embedded if viewing ${\mathcal X}$ as a subspace of ${\mathcal X}^{**}$ there is a projection $P$ of ${\mathcal X}^{**}$ onto ${\mathcal X}$ so that $\|z\| = \|Pz\| + \|z-Pz\|$ for all $z \in {\mathcal X}^{**}$ (see \cite{10}).

\bigskip
\noindent
{\bf 4.1. Proposition.} {\em 
Assuming $\Phi$ and $\Phi^*$ are mononorming and $k_{\Phi}(\tau) = 0$, ${\mathcal E}_*(\tau;\Phi)$ endowed with the norm $\|\ \|_M$ is $L$-embedded and separable.
}

\bigskip
\noindent
{\bf {\em Proof.}} Given $\varphi \in {\mathcal E}^*(\tau;\Phi)$ let $\varphi = \varphi_a + \varphi_s$ be its Lebesgue decomposition (Prop.~3.2) so that $\|\varphi\|_M = \|\varphi_a\| + \|\varphi_s\|_M$ and by the discussion in Subsection~$4.1$ we have $\varphi_a \in {\mathcal E}_*(\tau;\varphi)$. The fact that the map $\varphi \to \varphi_a$ is linear is a consequence of the uniqueness of the decomposition. Also clear is that $(\varphi_a)_a = \varphi_a$ so that defining $P\varphi = \varphi_a$ we have $P^2 = P$, $\|P\| \le 1$ and $\|\varphi|_M = \|\varphi_a\|_M + \|\varphi_s\|_M = \|P\varphi\|_M + \|(I-P)\varphi\|_M$. Since $\Phi^*$ is mononorming, ${\mathcal C}_1 \times ({\mathcal G}_{\Phi^*})^n$ is separable and so is then its quotient space ${\mathcal E}_*(\tau;\Phi)$.\qed

\subsection{Uniqueness}
\label{subsec4.3}

A theorem of H.~Pfitzner (\cite{12}) gives that an $L$-embedded Banach space has property $(X)$ of G.~Godefroy--M.~Talagrand, which in turn by the work of these authors \cite{7} implies strong uniqueness of that Banach space as a predual of its dual. Endowing the Banach space with an equivalent norm and its dual with the corresponding norm, any other predual of the dual will be isometric to the predual of the Banach space with the equivalent norm.

\bigskip
\noindent
{\bf 4.1. Corollary.} {\em 
Assume $k_{\Phi}(\tau) = 0$ and $\Phi,\Phi^*$ are mononorming. Then $({\mathcal E}_*(\tau;\Phi),\|\ \|)$ is the isometrically unique predual of $({\mathcal E}(\tau;\Phi),\||\ \||)$ and $({\mathcal E}_*(\tau;\Phi),\|\ \|_M)$ is the isometrically unique predual of $({\mathcal E}(\tau;\Phi),\||\ \||_M)$. More generally ${\mathcal E}(\tau;\Phi)$ has strongly unique predual.
}

\end{document}